\documentclass{article}

\usepackage[backend=bibtex, 
sortcites=true,
firstinits=true,
hyperref,
maxbibnames=99,
]{biblatex}

\bibliography{Bibliography}{}

\usepackage{enumerate}

\usepackage[centertags]{amsmath}
\usepackage{amsfonts}
\usepackage{amssymb}
\usepackage{amsthm}
\usepackage{newlfont}
\usepackage{mathtools}
\usepackage[top=1in, bottom=1in, left=1in, right=1in]{geometry}
\usepackage{tikz}

\theoremstyle{plain}
\newtheorem{theorem}{Theorem}[section]

\newtheorem{lemma}[theorem]{Lemma}

\theoremstyle{definition}
\newtheorem{definition}[theorem]{Definition}

\newcommand{\R}{\mathbb{R}}
\newcommand{\N}{\mathbb{N}}

\newcommand{\Deg}{\operatorname{Deg}}

\newcommand{\diam}{\operatorname{diam}}

\newcommand{\diameff}{\diam_{\operatorname{eff}}}

\newcommand{\Hidden}[1]{}

\begin{document}

\title{Reflective Graphs, Ollivier curvature, effective diameter, and rigidity}
\author{

Florentin M\"unch\footnote{MPI MiS Leipzig, muench@mis.mpg.de}
}
\date{\today}
\maketitle

\begin{abstract}
We give a discrete Bonnet Myers type theorem for the effective diameter assuming positive Ollivier curvature. We prove that this diameter bound is attained if and only if the graph is
a cocktail party graph, a Johnson graph, a halved cube, a Schläfli graph, a Gosset graph, or a cartesian product of the mentioned graphs with same Ollivier curvature.

As a key step in the proof, we introduce the notion of reflective graphs as graphs such that for any two neighbors there exists a certain self-inverse automorphism mapping one neighbor to another. We classify these graphs as arbitrary cartesian products of the graphs mentioned before. 
\end{abstract}


\section{Introduction}

Let $G=(V,E)$ be a simple, finite, connected graph. We define the effective diameter
\[
\diameff(G) = \frac 1 {|V|^2} \sum_{x,y \in V} d(x,y)
\]
where $d$ is the combinatorial graph distance.
We always assume that graphs are simple, finite and connected.
The Ollivier Ricci curvature $\kappa$ of an edge $x\sim y$ is given by (see \cite{munch2017ollivier,ollivier2007ricci,lin2011ricci})
\[
\kappa(x,y) := \inf_{\substack{f(y)-f(x)=1\\ \|\nabla f\|_\infty = 1}} \Delta f(x) - \Delta f(y)
\]
where $\|\nabla f\|_\infty:= \max_{x\sim y} |f(y)-f(x)|$ and $\Delta: \R^V \to \R^V$,
\[
\Delta f(x):=\sum_{y\sim x} (f(y)- f(x)).
\]
It is shown in \cite{jost2021characterizations} that Ollivier and Forman curvature coincide, a curvature notion on CW complexes defined via a discrete Bochner formula.

The degree $\Deg(x)$ of a vertex $x \in V$ is defined as $\Deg(x):= |\{y \in V: y \sim x\}|$.
It is well known that a positive lower curvature bound $\kappa \geq K$ implies an upper diameter bound, see \cite{ollivier2009ricci,liu2016bakry,munch2017ollivier}, given by
$
\diam(G):=\max_{x,y}d(x,y) \leq \frac{2\max \Deg}{K}.
$
In \cite{cushing2018rigidity} it was investigated for which graphs  equality holds. It turned out that under the additional assumption that for all $x\in V$ there exists $y \in V$ with $d(x,y)=\diam(G)$, the graphs for which holds equality are precisely the following:
 Cocktail party graphs,
 Johnson graphs $J(2n,n)$,
 halved cubes on $2\cdot 4^n$ vertices,
 Gosset graph, and
 Cartesian products of the mentioned graphs with same curvature.
It is still open whether the additional assumption can be dropped, and this article might be useful for solving this question.
Specifically, in this paper, we give a diameter bound for the effective diameter and classify all graphs for which this bound is attained. 
We now give the diameter bound.

\begin{theorem}
Let $G=(V,E)$ be a graph with $\kappa(x,y) \geq K >0$ for all $x\sim y$. Then,
\[
\diameff(G) \leq \frac{\max \Deg}{K}.
\]
\end{theorem}
This theorem reappears as Theorem~\ref{thm:EffDiamBound}.
We say a graph is effective Bonnet Myers sharp if the effective diameter bound holds with equality. 
We now characterize the effective Bonnet Myers sharp graphs.
Throughout the paper, we use the abbreviation "T.f.a.e." for "The following are equivalent".
\begin{theorem}
Let $G=(V,E)$ be a graph with $0<K:=\min_{x\sim y} \kappa(x,y)$. T.f.a.e.:
\begin{enumerate}[(i)]
\item $\diameff(G) = \frac{\max \Deg}{K}$,
\item $G$ is a graph from the following list:
\begin{enumerate}[(1)]
\item Cocktail party graph,
\item Johnson graph, 
\item Halved cube,
\item Schläfli graph,
\item Gosset graph,
\item Cartesian products of the above graphs with same curvature.
\end{enumerate}
\end{enumerate}
\end{theorem}
This theorem reappears as Theorem~\ref{thm:effBMsharpChar}.
It turns out that the effective Bonnet Myers sharp graphs all share a global combinatorial property which we call reflectiveness and discuss in the next section. We will show that reflectiveness almost characterizes effective Bonnet Myers sharpness. The only difference is that reflectiveness is compatible with arbitrary cartesian products, but effective Bonnet Myers sharpness is only compatible with cartesian products when the factors have the same curvature.
The key argument for classifying reflective graphs is the following,
\begin{itemize}
\item Locally disconnected $\Rightarrow$ Cartesian product
\item Locally connected $\Rightarrow$ Distance transitive 
\end{itemize}
Here, locally connected means that the induced subgraph of the neighborhood of any vertices is connected.
We then show that in the locally connected case, the graph has to be Lichnerowicz sharp, i.e. $K=\lambda$ where $K$ is the minimal curvature, and $\lambda$ is the smallest positive eigenvalue of $-\Delta$. We then use the classification of distance regular Lichnerowicz sharp graphs with an additional spectral condition from \cite[Theorem~6.5]{cushing2018rigidity} to classify reflective graphs.

In the appendix, we briefly explain what happens when replacing Ollivier by Bakry Emery curvature. It turns out that this is quite easy. Particularly, one gets the same diameter bound as in the Ollivier case. For Bakry Emery curvature however, this bound is attained only for hypercubes, and this almost immediately follows from \cite{liu2017rigidity}.

\subsection{Historical background}

Curvature is one of the fundamental concepts in geometry. 
Curvature bounds have many important analytic, geometric and topological implications. For example, the Bonnet Myers diameter bound states that a Riemannian manifold with a positive lower Ricci curvature bound can have at most the diameter as the round sphere with the same curvature \cite{myers1941riemannian}. Moreover, Cheng's sphere theorem states that this diameter bound is only attained for the round sphere \cite{cheng1975eigenvalue}.

In the last decade, there has been increasing interest in discrete Ricci curvature notions. There are four different basic discrete Ricci curvature notions, Ollivier curvature \cite{ollivier2007ricci,ollivier2009ricci}, Bakry Emery curvature \cite{schmuckenschlager1998curvature,lin2010ricci}, Forman curvature \cite{forman2003bochner,jost2021characterizations}, and entropic curvature \cite{erbar2012ricci,mielke2013geodesic}, and all come together with several sometimes non-linear modifications \cite{bauer2015li,munch2014li,dier2017discrete,
kempton2020large,eidi2020ollivier,ikeda2021coarse,
asoodeh2018curvature,devriendt2022discrete,
topping2021understanding}. 

Due to the lack of a chain rule for the discrete Laplacian, these curvature notions are all different. In this paper, we focus on the Ollivier curvature. Most of the curvature notions allow for a discrete Bonnet Myers theorem, i.e., an upper diameter bound under a positive Ricci bound, and this diameter bound is sharp for Ollivier and Bakry Emery.
In the attempt to classify the graphs for which the diameter bound is attained, interesting relations to spectral graph theory were discovered \cite{cushing2018rigidity}.
Particularly, it is conjectured that all these graphs are strongly spherical.
While strongly spherical graphs are fully classified \cite{koolen2004structure}, the classification of spherical graphs is still open.

It turns out that this paper is right at the interface between these questions:

We give a complete classification of all graphs for which the effective diameter bound is attained. To do so, we introduce the notion of reflective graphs. 
It turns out that the reflective graphs are a strict subclass of spherical graphs and a strict superclass of strongly spherical graphs.

\subsection{Structure of the paper}

The paper is structured as follows.
In the next section, we introduce and classify reflective graphs. In Section~\ref{sec:Curvature}, we show that effective Bonnet Myers sharp graphs are reflective. Finally in Section~\ref{sec:Classifications}, we summarize the results and characterize effective Bonnet Myers sharp graphs, reflective graphs, and distance regular Lichnerowicz sharp graphs, improving Theorem~6.5 in \cite{cushing2018rigidity}.
In the appendix, we briefly give effective diameter bounds in terms of Bakry Emery curvature and characterize equality.
\section{Reflective Graphs}\label{sec:Reflective}

Let $G=(V,E)$ be a simple, finite graph. Let $x\sim y$ be vertices. We write
\begin{itemize}
\item $V_x^y := \{x':d(x',x)<d(x',y) \}$,
\item $V^{xy} := \{z:d(z,x)=d(z,y) \}$.
\end{itemize}
We remark that $(V_x^y,V_y^x,V^{xy})$ is a partition of $V$.
\begin{definition}
A map $\phi : V\to V$ is called \emph{reflection} from $x$ to $y$ if
\begin{enumerate}[(i)]
\item $\phi$ is a graph automorphism,
\item $\phi^2 = id$,
\item $\phi(x)=y$,
\item $E(V_x^y,V_y^x) = \{(x',\phi(x')) : x' \in V_x^y\}$,
\item $\phi(z)=z$ for all $z \in V^{xy}$.  
\end{enumerate}
We say a graph is reflective, if for all $x\sim y$, there exists a reflection from $x$ to $y$.
We remark that this reflection from $x$ to $y$ is unique, and we call it $\phi_{xy}$.
\end{definition}

\begin{definition}
For vertices $x\sim y$,  and $x'\sim y'$, we write $(x,y) \parallel (x',y')$ if $x' \in V_x^y$ and $y' \in V_y^x$. In this case, we say the edges are parallel. 
\end{definition}

We remark that this relation depends on the  order of the pairs. It is easy to check that $(x,y) \parallel (x',y')$ if and only if $x' \in V_{x}^y$ and $y' = \phi_{xy}(x')$.

The aim of this section is to classify reflective graphs.
In the next subsection, we give fundamental properties of parallel edges. Afterwards we show that locally disconnected reflective graphs are cartesian products, and that locally connected reflective graphs are distance transitive. Our next step is to prove that locally connected reflective graphs are Lichnerowicz sharp, and we then apply the classification of Lichnerowicz sharp distance regular graphs from \cite{cushing2018rigidity}.
 We remark that a similar dichotomy between locally connected and disconnected graphs also appears in the classification of strongly spherical graphs \cite{koolen2004structure}. 
\subsection{Parallel edges}

We first give a useful characterization of parallel edges.
\begin{lemma}\label{lem:ParallelEquivalence}
Let $G=(V,E)$ be reflective. Let $x\sim y$ and $x'\sim y'$ be vertices. T.f.a.e.:
\begin{enumerate}[(i)]
\item $(x,y) \parallel (x',y')$,
\item $V_x^y = V_{x'}^{y'}$ and $V_y^x = V_{y'}^{x'}$.
\end{enumerate}
Particularly, if any of the equivalent conditions is satisfied, then \[
d(x,z)-d(y,z) = d(x',z)-d(y',z) \quad \mbox{ for all } z \in V.
\]

\end{lemma}
\begin{proof}
We first show $(i) \Rightarrow (ii)$, and particularly, the first equation in $(ii)$.
Suppose $z \in V_{x}^{y}$. We aim to show $z \in V_{x'}^{y'}$.
Let $(x_k)_{k=0}^n$ be a shortest path from $x$ to $z$.
We observe that $x_k \in V_{x}^y$. Let $y_k := \phi_{xy}(x_k) \in V_y^x$.

We show via induction that $x_k \in V_{x'}^{y'}$ and $y_k \in V_{y'}^{x'}$. The induction base $k=0$ is clear. Now, assume $x_{k-1} \in  V_{x'}^{y'}$ and $y_{x-1} \in V_{y'}^{x'}$.

Suppose $x_k \in V^{x'y'}$ for which we find a contradiction. Then, $x_{k} \sim x_{k-1}$ implying
\[
x_k = \phi_{x'y'}(x_k) \sim \phi_{x'y'}(x_{k-1}) =y_{k-1}. 
\]
This is a contradiction to $(iv)$ in the definition of reflections, as $y_{k-1} \in V_y^x$ has two different neighbors in $V_x^y$, namely $x_k$ and $x_{k+1}$, and the contradiction shows $x_k \notin V^{x'y'}$.

Next, suppose $x_k \in V_{y'}^{x'}$ for which we find a contradiction. Then, $x_{k-1} \in V_{x'}^{y'}$ has two different neighbors in $V_{y'}^{x'}$, namely $x_k$ and $y_{k-1}$. We remark that they are different as $x_k \in V_x^y$ and $y_{k-1} \in V_y^x$. This is a contradiction  to $(iv)$ in the definition of reflections.

The above two cases show $x_k \in V_{x'}^{y'}$. Similarly, $y_k \in V_{y'}^{x'}$ finishing the induction. Thus, $z \in V_{x'}^{y'}$. This proves $V_x^y \subseteq V_{x'}^{y'}$. The reverse inclusion follows analogously, showing $V_x^y = V_{x'}^{y'}$. By interchanging the role of $x$ and $y$, the second equation in $(ii)$ follows immediately, 
proving $(i) \Rightarrow (ii)$.

We finally prove $(ii) \Rightarrow (i)$.
By $(ii)$, we have $x' \in V_{x'}^{y'}=V_x^y$ and $y' \in V_{y'}^{x'}=V_{y}^x$ proving $(ii) \Rightarrow (i)$ and finishing the proof of the equivalence.
The "particularly" statement immediately follows from $(ii)$.
\end{proof}

We next show that one can find parallel edges close to every vertex. We write $B_1(x)=\{y :d(x,y) \leq 1\}$ for the balls of radius one.
\begin{lemma}\label{lem:ParellelToBall}
Let $x\sim y \in V$ and $z \in V$. Then, there exist $x',y' \in B_1(z)$ such that
\[
(x,y) \parallel (x',y').
\]
\end{lemma}

\begin{proof}
Suppose not. 
If $z \in V_x^y$, then $(x,y) \parallel (z,\phi(z))$.
If $z \in V_y^x$, then $(x,y) \parallel (\phi(z),z)$.
Hence, we can assume  $z \in V^{xy}$.
By induction, we assume $n:=d(x,z)$ to be minimal.
Let $p \sim z$ on a geodesic from $z$ to $x$, i.e. $d(p,x)=n-1$.
We notice
\[
n-1 \leq d(p,y) \leq n.
\]
\begin{description}
\item[Case 1: $d(p,y)=n-1$.]

Let $x' := \phi_{pz}(x)$, and $y' :=\phi_{xy}(x')$. 
Then $n-1=d(x',z) = d(y',z)$. Moreover $(x',y') \parallel (x,y)$ contradicting minimality of $n$ as "$\parallel$" is an equivalence relation by Lemma~\ref{lem:ParallelEquivalence}.

\item[Case 2: $d(p,y)=n$.]

Let $y' :=\phi_{xy}(p)$. As $p \sim z \in Z$, we have $y' \sim z$. Hence, $(x,y) \parallel (p,y')$ and $p,y' \in B_1(z)$ being a contradiction.
\end{description}
Now the proof is finished by contradiction.
\end{proof}

\subsection{Cartesian products}

We aim to show that locally disconnected graphs are cartesian products. We recall that locally disconnected means that
the induced subgraph on the neighborhood of a vertex is disconnected. 
We now give a characterization of cartesian product graphs from \cite{MathoverflowCartesianProducts,imrich2007recognizing}.
\begin{theorem}\cite[Theorem~3.1]{imrich2007recognizing} \label{thm:CartesianFacorization}
Let $G=(V,E)$ be a graph. T.f.a.e.:

\begin{enumerate}[(i)]
\item $G$ is a cartesian product of two non-trivial graphs,
\item The graph $(E,\sim)$ is disconnected where $(x,y) \sim (x',y')$ iff
\begin{enumerate}[(a)]
\item $d(x,x')+d(y,y')\neq d(x,y')+ d(y,x')$ or
\item $\{x\}=\{x'\}= S_1(y) \cap S_1(y')$ maybe after interchanging $x$ and $y$ or $x'$ and $y'$.
\end{enumerate}
\end{enumerate}
\end{theorem}
We now investigate the connection between the relation "$\sim$" from the theorem above and the parallel relation "$\parallel$" for reflective graphs.
\begin{lemma}\label{lem:parallelSim}
Let $G=(V,E)$ be a reflective graph.
Let $(x,y) \parallel (x',y') \in E$ and $(v,w) \in E$. T.f.a.e.:
\begin{enumerate}[(i)]
\item $d(x,v)+d(y,w)\neq d(x,w)+d(y,v)$,
\item $(x,y) \sim (v,w)$,
\item $(x',y') \sim (v,w)$.
\end{enumerate}
\end{lemma}
\begin{proof}
For reflective graphs, condition (b) will never be decisive as $|S_1(y) \cap S_1(y')| \geq 2$ whenever $d(y,y')=2$, and if $d(y,y')=1$, then (a) is satisfied as well. 
This proves $(i)\Rightarrow (ii)$.
We now prove $(i) \Rightarrow (iii)$.
By Lemma~\ref{lem:ParallelEquivalence} for $z \in \{v,w\}$,
\[
d(z,x) - d(z,y) = d(z,x') - d(z,y')
\]
Therefore,
\begin{align*}
d(x',v) + d(y',w) &= d(v,x) - d(v,y) + d(v,y')  +  d(w,x') + d(w,y) - d(w,x) \\
&\neq  d(v,y')  +  d(w,x')
\end{align*}
proving $(i) \Rightarrow (iii)$. The reverse direction $(iii) \Rightarrow (i)$ can be proven similarly.
\end{proof}

We next show that locally disconnected graphs are cartesian products. We write $S_n(x)=\{y:d(x,y)=n\}$ for the spheres.
\begin{lemma}\label{lem:LocDisconnectedImpliesCartesianProduct}
Let $G=(V,E)$ be a reflective graph. Let $z \in V$. Suppose the induced subgraph of $S_1(z)$ is disconnected. Then, $G$ is a cartesian product with non-trivial factors.
\end{lemma}

\begin{proof}
Let $p,p' \in S_1(z)$ in different connected components.
Suppose, $G$ is not a cartesian product. Then, by Theorem~\ref{thm:CartesianFacorization}, there is a path $(z_k,p_k)_{k=0}^n$ with respect to $\sim$ and $z_0=z_n=z$ and $p_0=p$ and $p_n=p'$.
By Lemma~\ref{lem:ParellelToBall}, there exist $(z_k',p_k') \in B_1(z)$ with $(z_k',p_k') \parallel (z_k,p_k)$. By Lemma~\ref{lem:parallelSim}, this is also a path with respect to $\sim$.
Let $(x',y'):=(p_k',z_k')$ be the first edge with a vertex outside of $C \cup \{z\}$ where $C$ is the connected component of $p$ within $S_1(z)$. Then, $(x,y):=(z_{k-1}',p_{k-1}') \in C \cup \{z\}$.
We now show via case distinction that $(x,y) \not\sim (x',y')$ which will be a contradiction.

We will use that if $x\neq z \neq x'$, then, $d(x,x')=2$ as $x \in C$ but $x' \notin C$. Moreover, if $x',y' \neq z$, then both $x',y' \notin C$ as $x' \sim y'$.
\begin{description}
\item[Case 1:] $x=z=x'$. Then, 
$d(x,x') + d(y,y')=0+2=1+1=d(x,y')+d(y,x')$.

\item[Case 2:] $x=z$ and $x',y' \neq z$. Then, $d(x,x') + d(y,y')=1+2=1+2=d(x,y')+d(y,x')$.

\item[Case 3:] $z\neq x,y$ and $x'=z$. Then,
$d(x,x') + d(y,y')=1+2=2+1=d(x,y')+d(y,x')$.

\item[Case 4:] $z\neq x,y$ and $x',y' \neq z$. Then,
$d(x,x') + d(y,y')=2+2=2+2=d(x,y')+d(y,x')$.
\end{description}
The remaining cases are analogous.
This finishes the case distinction giving the contradiction finishing the proof. 
\end{proof}

We finally prove that a cartesian product is reflective if and only if all factors are reflective.

\begin{lemma}\label{lem:CartesianProductInheritsReflective}
Let $G_i=(V_i,E_1)$ be graphs for $i=1,2$. T.f.a.e.:
\begin{enumerate}[(i)]
\item $G_1$ and $G_2$ are reflective,
\item $G_1\times G_2$ is reflective.
\end{enumerate}
\end{lemma}
\begin{proof}
We first prove $(i)\Rightarrow (ii)$. Let $(x,y)=((x_1,x_2),(y_1,y_2))$ be an edge in $G_1 \times G_2$. W.l.o.g., $x_2=y_2$, as we can interchange $G_1$ and $G_2$. Then, $x_1 \sim y_1$, and we have a reflection $\phi_{x_1y_1}$ in $G_1$. We extend this to
\[
\phi(z_1,z_2) := (\phi_{x_1y_1}(z_1), z_2).
\]
It is straight forward to check that $\phi$ is a reflection. This proves $(i) \Rightarrow (ii)$.

We finally prove $(ii) \Rightarrow (i)$.
We show that $G_1$ is reflective. Let $x_1 \sim y_1 \in V_1$ and $x_2=y_2 \in V_2$.
Then, there is a reflection $\phi$ from $(x_1,x_2)$ to $(y_1,y_2)$. We notice
\[
V_{(x_1,x_2)}^{(y_1,y_2)} = (V_1)_{x_1}^{y_1} \times V_2.
\]
Hence, for all $z_1 \in V_1$, there exists $z_1' \in V_1$ such that for all $z_2 \in V_2$, 
\[
\phi(z_1,z_2) = (z_1',z_2).
\]
We define $\phi_{x_1y_1}(z):=z'$.
Again, it is straight forward to check that $\phi_{x_1y_1}$ is a reflection. This proves $(ii) \Rightarrow (i)$ and finishes the proof.
\end{proof}

\subsection{Locally connected implies distance transitive}
Having investigated the locally disconnected case, 
we now come to the locally connected case, and make use of the automorphisms $\phi_{xy}$ to show distance transitivity.
We first show that the class of reflective graphs is closed under taking induced subgraph on $V_x^y$ for $x\sim y$.
We recall that a subset  $W\subseteq V$ is convex if every shortest path with first and last vertex in $W$ must have all vertices in $W$.

\begin{lemma}\label{lem:VxyIsometricReflective}
Let $G=(V,E)$ be reflective and $x \sim y$.
Then, the induced subgraph on $V_x^y$ is convex and reflective.
\end{lemma}

\begin{proof}
We first show that $V_x^y$ is convex. 
Let $v=v_0 \sim \ldots \sim v_n = v'$ be a shortest path in $G$.
Let $w:= \phi_{xy}(v)$. By Lemma~\ref{lem:ParallelEquivalence}, we have $v' \in V_x^y=V_{v}^w$. Thus, $(w,v_0,\ldots,v_n)$ is a shortest path showing $v_k \in V_v^w = V_x^y$ showing that every shortest path from $v$ to $v'$ stays in $V_x^y$.

We now show that the induced subgraph on $V_x^y$ is reflective. We first prove that $\phi_{vw} : V_x^y \to V_x^y$ for all neighbors $v,w \in V_x^y$.
Let $v' \in V_x^y$. We want to show $\phi_{vw}(v') \in V_x^y$
If $v' \in V^{vw}$, then, $\phi(v')=v' \in V_x^y$. 
If $v' \in V_w^v$, we interchange $v$ and $w$, so we are left with the case $v' \in V_v^w$. Let $w'=\phi_{vw}(v')$. Then, $(v,w) \parallel (v',w')$. By Lemma~\ref{lem:ParallelEquivalence},
\begin{align*}
d(v,x)-d(w,x) &= d(v',x)-d(w',x), \mbox{ and }\\
d(v,y)-d(w,y) &= d(v',y)-d(w',y).
\end{align*}
Subtracting the equations gives
$d(w',y)-d(w',x)=1$ as $v,w,v' \in V_x^y$. Thus, $w' \in V_x^y$.
Thus, $\phi_{vw}$ maps $V_x^y$ to itself as claimed.
As $V_x^y$ is an isometric subgraph, we see that 
$\left(V_x^y\right)_v^w = V_v^w \cap V_x^y$ and $\left(V_x^y\right)_w^v = V_w^v \cap V_x^y$. Thus,
$\phi_{vw}|_{V_x^y}$ is a reflection on $V_x^y$. This proves that the induced subgraph on $V_x^y$ is reflective as $v \sim w \in V_x^y$ were chosen arbitrarily. This finishes the proof.
\end{proof}

We now prove that the neighborhood of a vertex is isometric.
We recall that a subset $W\subseteq V$ is called isometric if for all $w,w'$, there exists a geodesic from $w$ to $w'$ entirely in $W$. 
\begin{lemma}\label{lem:S1isometric}
Let $G=(V,E)$ reflective and $x \in V$. Assume $S_1(x)$ is connected. Then, $S_1(x)$ is isometric.
\end{lemma}
\begin{proof}
Suppose not. Then, there exist vertices in $S_1(x)$ such that a shortest path within $S_1(x)$ between them has length at least three. Let the first four vertices of such a path be denoted by $y\sim z \sim z' \sim y' \in S_1(x)$.
We aim to find $w\in S_1(x)$ with $y \sim w \sim y'$ which would be a contradiction.
We define
\[
v:=\phi_{xy}(z') = \phi_{xz'}(y).
\]
As $d(x,z')=2$, we get
$v \sim y,z,z'$ and $v \not\sim x$. 
Moreover, $v \not \sim y'$ as otherwise, $v \in V_y^x$ would have two neighbors in $V_x^y$, namely $y'$ and $z'$.
Similar to $v$, we define
\[
v':=\phi_{xy'}(z) = \phi_{xz}(y').
\]
By similar arguments as before, we get $v' \sim z,z',y'$ and $z \not\sim x,y$.
As $d(x,v)=d(y',v) = 2$, and as $v\sim z$, we get
\[
v= \phi_{xy'}(v) \sim \phi_{xy'}(z) = v'.
\]
We finally define 
\[
w:= \phi_{xz}(v).
\]
Obviously, $w \in S_1(x)$. As $v\sim y$, we get
\[
y=\phi_{xz}(y) \sim \phi_{xz}(v) = w.
\] 
As $v\sim v'$, we get
\[
w = \phi_{xz}(v) \sim \phi_{xz} (v') = y'
\]
as desired. This contradicts the assumption that shortest paths within $S_1(x)$ between $y$ and $y'$ have length three.
The contradiction finishes the proof.
\end{proof}

We now show that intersections of the neighborhood of a vertex with certain spheres are isometric, and particularly, connected. 
\begin{lemma}\label{lem:connectedIntersection}
Let $G=(V,E)$ be reflective. Let $x,x' \in V$ with $d(x,x')=n-1$ for some $n\geq 1$. Assume $S_1(x)$ is connected. Then,
\[
S_1(x) \cap S_{n}(x') \mbox{ is isometric.}
\]
\end{lemma}

\begin{proof}
Suppose not. By induction, we can assume $n$ to be minimal. By Lemma~\ref{lem:S1isometric}, we can assume $n \geq 2$.

Let $x''\sim x'$ be on a shortest path from  $x'$ to $x$, i.e. $d(x,x'')=n-2$.
We notice
\[
S_1(x) \cap S_{n}(x') = S_1(x) \cap S_{n-1}(x'') \cap V_{x''}^{x'}. 
\]
By induction hypothesis, we see that $S_1(x) \cap S_{n-1}(x'')$ is isometric. As $V_{x''}^{x'}$ is convex by Lemma~\ref{lem:VxyIsometricReflective}, it follows that  
$S_1(x) \cap S_{n-1}(x'') \cap V_{x''}^{x'}$ is also isometric as an intersection of an isometric and a convex set.
This finishes the induction step and thus, the proof is complete.
\end{proof}

We finally show distance transitivity by suitably composing the automorphisms $\phi_{xy}$.
\begin{lemma}\label{lem:DistTransitive}
Let $G=(V,E)$ is reflective and locally connected. Then, $G$ is distance transitive.
\end{lemma}

\begin{proof}
We first show vertex transitivity.
Let $x,x' \in V $ and let $x=x_0\sim\ldots = x_n = x'$ be a path. Then, the automorphism $\phi_{x_{n}x_{n-1}} \circ \ldots \circ \phi_{x_1 x_0}$ maps $x$ to $x'$ showing vertex transitivity.

We now show distance transitivity.
Let $n \in \N$.
Let $x,y,x',y' \in V$ with $d(x,y)=d(x',y')=n$. We aim to find an automorphism $\psi$ mapping $x$ to $x'$ and $y$ to $y'$. We proceed via induction over $n$. The case $n=0$ is clear by vertex transitivity.

Now for $n>0$, we assume the claim to be true for $n-1$.
Let $z\sim x$ on a shortest path from $x$ to $y$ and $z'\sim y'$ on a shortest path from $x'$ to $y'$. By induction assumption, we can assume $x=x'$ and $z=z'$. Hence,
$y,y' \in S_1(z) \cap S_n(x)$. By Lemma~\ref{lem:connectedIntersection}, there is a path $y=y_0 \sim \ldots \sim y_m=y'$ in $S_n(x)$.
Now we define
\[
\psi:=\phi_{y_{m}y_{m-1}}\circ \ldots \circ \phi_{y_1y_0}.
\]
Clearly, $\psi(y) = \psi(y')$. Moreover as $d(y_k,x)=n$, we have
$\phi_{y_ky_{k-1}}(x)=x$ for all $k$ showing that $\psi(x)=x=x'$.
This finishes the induction, and the proof is now complete.
\end{proof}

\subsection{Lichnerowicz sharpness}
In order to classify reflective graphs, we want to use the classification of Lichnerowicz sharp, distance regular graphs from \cite[Theorem~6.5]{cushing2018rigidity}. We now give the details.
We recall
\[
\Delta f(x) = \sum_{y\sim x} (f(y)-f(x))
\]
and
$Lip(1):=\{f \in \R^V : f(y)-f(x) d(x,y) \mbox{ for all } x,y\in V\}$
and the Ollivier curvature $\kappa$ of an edge $x\sim y$ is
\[
\kappa(x,y) = \inf_{\substack{f \in Lip(1)\\f(y)-f(x)=1}} \Delta f(x)-\Delta f(y).
\]
The Lichnerowicz estimate states that 
\[
\lambda \geq \min \kappa
\]
where $\lambda$ is the smallest positive eigenvalue of $-\Delta$.
We show that reflective graphs are Lichnerowicz sharp, i.e., $K=\lambda$. For convenience, we restrict ourselves to locally connected graphs.

We recall that a graph $G$ is distance regular with intersection array $(b_0,b_1\ldots,b_{L-1}; 
c_1=1\ldots c_L)$ if it has diameter $L$, and if every vertex in $S_n(x)$ has $b_n$ neighbors in $S_{n+1}(x)$ and $c_n$ neighbors in $S_{n-1}(x)$, for arbitrary $x \in V$.
We now show that locally connected reflective graphs are Lichnerowicz sharp.

\begin{theorem}\label{thm:ReflectiveImpliesLichSharp}
Let $G=(V,E)$ be locally connected and reflective. 
Assume $G$ has intersection array 
\[(b_0,b_1\ldots,b_{L-1}; 
c_1=1\ldots c_L).\] Then, the non-normalized Ollivier curvature of all edges is
\[
\kappa = 1+ b_0 - b_1.
\]
Moreover, $G$ is Lichnerowicz sharp, i.e. the smallest positive eigenvalue $\lambda$ of the non-normalized Laplacian $-\Delta$ satisfies
\[
\lambda=\kappa.
\]

\end{theorem}
\begin{proof}
By the Lichnerowicz estimate $\lambda \geq \kappa$, it suffices to prove
\[\kappa \geq 1+ b_0 - b_1 \geq \lambda.\]
For showing $\lambda \leq 1+ b_0 - b_1$, we construct an eigenfunction $f$  with $\Delta f=-(1+b_0-b_1) f$.
Let $x\sim y$. Let
\[
f:=d(x,\cdot) - d(y,\cdot).
\]
We notice 
\[
b_1 = |S_2(x)\setminus{S_1(y)}| = |V_y^x \cap S_1(y)| = |V_x^y \cap S_1(x)|
\]
and
\[
b_0 = 1 + |V_x^y \cap S_1(x)| + |V^{xy} \cap S_1(x)|.
\]
Thus,
\begin{align*}
\Delta f(x) &= f(y)-f(x) + \sum_{z \in V^{xy}\cap S_1(x)}(f(z) - f(x)) + \sum_{z \in V_x^y\cap S_1(x)}(f(z)-f(x)) \\&= 2 + |V^{xy} \cap S_1(x)| 
\\&= 1+ b_0 - b_1 = -f(x)(1+ b_0 - b_1).
\end{align*}
Similarly, $\Delta f(y) = - f(y)(1+ b_0 - b_1)$.
For $x' \in V_x^y$, we set $y'=\phi_{xy}(x')$ and get $f=d(x',\cdot)-f(y',\cdot)$ by Lemma~\ref{lem:ParallelEquivalence} showing $\Delta f(x') = -f(x')(1+ b_0 - b_1)$, and similarly for $y' \in V_y^x$.
For $z \in V^{xy}$, we have with $\phi=\phi_{xy}$,
\[
\Delta f(z) = \Delta f\circ \phi (z) = \Delta (d(y,\cdot)-\Delta(x,\cdot))(z) = - \Delta f(z)
\]
giving $\Delta f(z)=0=f(z)$. This shows that $f$ is an eigenfunction proving $1+b_0-b_1 \geq \lambda$.

We finally prove $\kappa \geq 1 + b_0 - b_1$.
By distance transitivity, the curvature is constant. Let $x\sim y$.
Let $f \in Lip(1)$ with $f(y)-f(x)=1$. We aim to show
\[
\Delta f(x)-\Delta f(y) \geq 1+b_0 - b_1.
\]
Let $\phi:=\phi_{xy}$.
We notice $\Delta f(y) =  \sum_{z\sim x} f(\phi(z))-f(y)$.
Hence,
\[
\Delta f(x)-\Delta f(y) = (f(y)-f(x))b_0 + \sum_{z \sim x} f(z)- f(\phi(z)). 
\]
We observe
\begin{align*}
\sum_{z\sim x} f(z)- f(\phi(z))  = f(y)-f(x) + \sum_{x\sim z \in V_x^y } f(z) - f(\phi(z)) + \sum_{x\sim z \in V^{xy}} f(z)-f(\phi(z)) .
\end{align*}
The latter sum vanishes as $z=\phi(z)$ on $V^{xy}$. The second sum can be bounded from below by $-|S_1(x) \cap V_x^y|=b_1$ as $z\sim \phi(z)$ and as $f \in Lip(1)$. Thus,
\[
\sum_{z\sim x} f(z)- f(\phi(z))   \geq 1 - b_1
\]
giving
\[
\Delta f(x)-\Delta f(y) \geq 1+b_0 - b_1
\]
where we used $f(y)-f(x)=1$. This proves $\kappa \geq 1+b_0-b_1$ and finishes the proof of the theorem.
\end{proof}

\subsection{First step to classification}
In this subsection, we aim to show that reflective graphs are cartesian products of Cocktail party graphs, Johnson graphs, halved cubes, the Schläfli graph, and the Gosset graph. We will prove in Section~\ref{sec:Classifications} that all graphs from the list are indeed reflective.
We will use the characterization of distance regular Lichnerowicz sharp graphs with second largest adjacency matrix eigenvalue $\theta = b_1 -1$ from \cite{cushing2018rigidity} which we recall now.
\begin{theorem}[Theorem~6.5 in \cite{cushing2018rigidity}]\label{thm:LichSharpChar}
Let $G=(V,E)$ be distance regular with intersection array
\[
(b_0,b_1,\ldots,b_{L-1};1=c_1,\ldots,c_L)
\]
Assume the second largest eigenvalue $\theta$ of the adjacency matrix satisfies $\theta = b_1 - 1$. Moreover assume $G$ is Lichnerowicz sharp, i.e., $\min \kappa = \lambda$ where $\lambda$ is the smallest positive eigenvalue of $-\Delta$. Then, $G$ is a graph of the following list:
\begin{enumerate}[(1)]
\item Cocktailparty graphs $CP(k)$,
\item Hamming graphs $K_n^m$,
\item Johnson graphs $J(n,k)$,
\item Halved cubes $Q^{(2)}(n)$,
\item Schlaefli graph,
\item Gosset graph.
\end{enumerate}
Conversely, all graphs from the list satisfy the properties above.
\end{theorem}
We remark that in \cite{cushing2018rigidity}, Lichnerowicz sharpness refers to the normalized Laplacian which due to constant vertex degree is equivalent to our version.
Moreover, we note that Hamming graphs are cartesian products of Johnson graphs. We now apply the above theorem to reflective graphs.

\begin{lemma}\label{lem:locConnectedAndReflectiveImpliesList}
Let $G=(V,E)$ be locally connected and reflective. Then, the graph $G$ is a graph from the list in Theorem~\ref{thm:LichSharpChar}.
\end{lemma}
\begin{proof}
By Lemma~\ref{lem:DistTransitive}, the graph is distance regular. By Theorem~\ref{thm:ReflectiveImpliesLichSharp}, the graph is Lichnerowicz sharp. As $\theta = b_0 - \lambda$, Theorem~\ref{thm:ReflectiveImpliesLichSharp} furthermore implies $\theta = b_1 - 1$. Hence, Theorem~\ref{thm:LichSharpChar} is applicable finishing the proof.
\end{proof}

As locally disconnected reflective graphs are cartesian products by Lemma~\ref{lem:LocDisconnectedImpliesCartesianProduct}, we can now classify all reflective graphs.
\begin{theorem}\label{thm:ReflectiveImpliesList}
Let $G=(V,E)$ be reflective. Then, the graph $G$ is a cartesian product of graphs in the list of Theorem~\ref{thm:LichSharpChar}.
\end{theorem}

\begin{proof}
Lemma~\ref{lem:LocDisconnectedImpliesCartesianProduct} and Lemma~\ref{lem:CartesianProductInheritsReflective} imply that $G$ is a cartesian product of locally connected, reflective graphs.
Lemma~\ref{lem:locConnectedAndReflectiveImpliesList} implies that all factors are graphs from the list. This finishes the proof.
\end{proof}


We remark that the class of reflective graphs lies strictly between the spherical and strongly spherical graphs investigated in \cite{koolen2004structure}. While strongly spherical graphs are fully classified, the classification of spherical graphs is still open, and this article might be useful for answering this question.

\section{Ollivier curvature and effective diameter}\label{sec:Curvature}

In this section, we give the effective diameter bound mentioned in the introduction, and show that equality implies that the graph is reflective. We first recall the effective diameter.
Let $G=(V,E)$ be a graph. Then, the effective diameter is
\[
\diameff=\diameff(G) = \frac 1 {|V|^2} \sum_{x,y \in V} d(x,y).
\]

\subsection{Diameter bounds}
In this subsection we prove the effective diameter bound using 
the long range Ollivier curvature
\[
\kappa(x,y) = \inf_{\substack{\|\nabla f\|_\infty = 1\\f(y)-f(x)=d(x,y)}} \frac{\Delta f(x)-\Delta f(y)}{d(x,y)}
\]
and using the fact that $\kappa(x,y) \geq K$ for all $x,y \in V$ as soon as this estimate holds for all neighbors $x\sim y$.
\begin{theorem}\label{thm:EffDiamBound}
Let $G=(V,E)$ be a graph with $\kappa(x,y) \geq K>0$ for all $x\sim y$. Then,
\[
\diameff \leq \frac{\max \Deg}{K}.
\]
\end{theorem}

\begin{proof}
We prove a stronger statement, namely that for all $x\in V$,
\[
\frac 1 {|V|}\max_x \sum_y d(x,y) \leq \frac{\max \Deg}{K}.
\]
Let $x\in V$ and $f:=d(x,\cdot)$.
Due to the curvature bound $\kappa \geq K$, we have for all $y \in V$,
\[
\Delta f(x) - \Delta f(y) \geq Kd(x,y).
\]
Summing up over all $y \in V$ gives
\[
|V|\Deg(x) = \sum_{y}(\Delta f(x)-\Delta f(y)) \geq K \sum_y d(x,y).
\]
Rearranging finishes the proof.
\end{proof}

\subsection{Rigidity and distance function as eigenfunction}

We now consider effective Bonnet Myers sharp graphs, i.e., graphs attaining the effective diameter bound, and aim to show that they are reflective.
This section is inspired by \cite[Section~3.3.2 and 3.3.3]{munch2019discrete}.
We first prove that $d(x,\cdot)$ is a shifted eigenfunction.
\begin{lemma}\label{lem:BMsharpDistanceEigenfunction}
Let $G=(V,E)$ be a graph with $\kappa(x,y)\geq K>0$ for all $x\sim y$. 
Suppose $\diameff = \frac {\max \Deg}K$.
Let $x \in V$. Then,
\[
\Delta d(x,\cdot) = \Deg(x)- Kd(x,\cdot).
\]
Moreover, $\Deg = \max \Deg$ and $\kappa=K$.
\end{lemma}

\begin{proof}
Let $x,y \in V$ and $f:=d(x,\cdot)$.
By the proof of Theorem~\ref{thm:EffDiamBound}, and due to equality in the effective diameter bound, we get
\[
\Delta f(x)-\Delta f(y) = K (f(y) - f(y)). 
\]
Choosing $x\sim y$ and using $\kappa \geq K$ yields $\kappa(x,y)=K$. As in the proof of Theorem~\ref{thm:EffDiamBound}, we estimated $\max \Deg \geq \Deg(x)$, we also get $\Deg(x)=\max \Deg$. We finally observe
\[
\Delta f(y) = \Delta f(x) - K(f(y)-f(x)) = \Deg(x) - Kd(x,y).
\]
This finishes the proof as $f=d(x,\cdot)$.
\end{proof}

\subsection{Rigidity and reflectiveness} We aim to show that effective Bonnet Myers sharp graphs are reflective.
We first give combinatorial implications of a positive curvature bound, namely that a positively curved edge has to be in many triangles, and there must be a perfect matching between the remaining neighbors. The following lemma is given in \cite[Proposition~3.1]{munch2019discrete}. For a version with the normalized graph Laplacian, also see \cite[Proposition~2.7]{cushing2018rigidity}. We recall $B_1(x)=\{y:d(x,y)\leq 1\}$.

\begin{lemma}[{{\cite[Proposition~3.1]{munch2019discrete}}}]\label{lem:curvTriangleMatch} 
Let $G=(V,E)$ be a graph. Let $x\sim y$. Then, $(x,y)$ contains at least $\kappa(x,y)-2$ triangles. In case of equality, there is a perfect matching between the remaining neighbors of $x$ and the remaining neighbors of $y$, i.e., there is a bijecitve map $\phi:B_1(x) \setminus B_1(y) \to B_1(y) \setminus B_1(x)$ such that $z \sim \phi(z)$ for all $z$ in the domain of $\phi$.
\end{lemma}

We next show that there is a perfect matching between $V_x^y$ and $V_y^x$. We recall $V_x^y$ is the set of vertices closer to $x$ than to $y$, and $V^{xy}$ is the set of vertices having same distance to $x$ and $y$.

\begin{lemma} \label{lem:BMsharpMatchingVxyVyx}
Let $G$ be effective Bonnet Myers sharp with curvature $K$. Let $x\sim y$. Then, for every $x' \in V_x^y$ there is exactly one $y' \in V_y^x$ with $x'\sim y'$, and exactly $K-2$ neighbors of $x'$ are in $V^{xy}$.
\end{lemma}

\begin{proof}
Let $D=\max \Deg$.
Let $x' \in V_x^y$. We notice $\Delta (d(x,\cdot) \wedge d(y,\cdot))(x) = D - 1$. Thus,
\[
\Delta (d(x,\cdot) \wedge d(y,\cdot))(x') \leq D - 1 - d(x,x')K
\]
where "$\wedge$" denotes the minimum.
By Lemma~\ref{lem:curvTriangleMatch}, the edge $(x,y)$ is contained in at least $K-2$ triangles, and thus,
\[
\Delta (d(x\cdot) \vee d(y,\cdot))(x) \leq D-1-(K-2) =  D+1-K
\]
implying
\[
\Delta (d(x,\cdot) \vee d(y,\cdot))(x') \leq D+1-K - d(x,x')K
\]
by the curvature bound, where "$\vee$" denotes the maximum.
By Lemma~\ref{lem:BMsharpDistanceEigenfunction},
\[
\Delta (d(x,\cdot) + d(y,\cdot))(x') = 2D - K(d(x,x')+d(y,x'))=2D-K -2Kd(x,x').
\]
Thus, adding up the above inequalities for $x'$, we get equality, meaning
\[
\Delta (d(x,\cdot) \wedge d(y,\cdot))(x') = D - 1 - d(x,x')K.
\]
Moreover by Lemma~\ref{lem:curvTriangleMatch}, we have
$\Delta d(x,\cdot)(x') = D-Kd(x,x')$, and thus,
\[
\Delta (1_{V_y^x})(x') = \Delta d(x,\cdot)(x') - \Delta (d(x,\cdot) \wedge d(y,\cdot))(x') = 1
\]
showing that $x'$ has exactly one neighbor in $V_y^x$.
Similarly, $\Delta d(y,\cdot)(x') = D-Kd(x,x')-K$ and thus,
\[
\Delta 1_{V_{x}^y}(x')=\Delta d(y,\cdot) (x') - \Delta (d(x,\cdot) \wedge d(y,\cdot))(x') = 1-K
\]
showing that $x'$ has exactly $K-1$ neighbors in $V\setminus V_x^y$, implying that $x'$ has exactly $K-2$ neighbors in $V^{xy}$.
This finishes the proof.
\end{proof}

We are now ready to give the main result of this section.
\begin{theorem}\label{thm:BMsharpImpliesReflective}
Every effective Bonnet Myers sharp graph is reflective.
\end{theorem}

\begin{proof}
Let $x\sim y$.
By the above theorem, there exists a perfect matching between $V_x^y$ and $V_y^x$ inducing a map $\phi$ mapping every $x' \in V_x^y$ to its unique neighbor in $V_y^x$, and conversely, and leaving all vertices in $V^{xy}$ fixed.
We aim to show that $\phi$ is a reflection.
All properties of a reflection are clear except that $\phi$ is an automorphism which we prove now.
Let $v\sim w$ be vertices. We aim to show that $\phi(v) \sim \phi(w)$. We proceed by case distinction.
\begin{description}
\item [Case 1:] $v \in V_x^y$ and $w \in V_y^x$. Then, $\phi(v)=w$ and $\phi(w)=v$, clearly showing  $\phi(v) \sim \phi(w)$.
\item [Case 2:] $v,w \in V^{xy}$. Then, $v= \phi(v)$ and $w=\phi(w)$, showing  $\phi(v) \sim \phi(w)$.
\item [Case 3:] $v \in V_x^y$ and $w \in V^{xy}$. By Lemma~\ref{lem:curvTriangleMatch}, the edge $(v,\phi(v))$ contains at least $K-2$ triangles. However, all common neighbors of $v$ and $\phi(v)$ must be in $V^{xy}$ by the first part of Lemma~\ref{lem:BMsharpMatchingVxyVyx}. By the second part of Lemma~\ref{lem:BMsharpMatchingVxyVyx}, $v$ has exactly $K-2$ neighbors in $V^{xy}$ showing that they must all be common neighbors of $v$ and $\phi(w)$. Specifically, $\phi(w)=w \sim \phi(v)$.
\item[Case 4:] $v,w \in V_x^y$. By the case above, the edge $(v,\phi(v))$ contains exactly $K-2$ triangles. By Lemma~\ref{lem:curvTriangleMatch}, there must be a perfect matching between the remaining neighbors of $v$ and $\phi(v)$. Specifically, there must be some $w' \notin B_1(v)$ with $w\sim w' \sim \phi(v)$. As $\phi(v)$ has at most one neighbor in $V_x^y$, namely $v$, infer $w' \notin V_x^y$. As all neighbors of $\phi(v)$ in $V^{xy}$ are also neighbors of $v$, but $w'$ is not a neighbor of $v$, we obtain $w' \in V_y^x$. As $\phi(w)$ is the unique neighbor of $w$ in $V_y^x$, we see that $w'=\phi(w)$. Particularly, $\phi(v) \sim \phi(w)$.
\end{description}
The remaining cases can be proven analogously. As $\phi=\phi^{-1}$, the case distinction shows that $\phi$ is a graph automorphism, and thus a reflection. As $x\sim y$ are chosen arbitrarily, the claim of the theorem follows immediately. 
\end{proof}

\subsection{Lichnerowicz sharpness}
We recall that Lichnerowicz sharp graphs are graphs with $\lambda = K$.
We also recall that distance regular Lichnerowicz sharp graphs with an additional spectral condition are precisely the graphs in the list in Theorem~\ref{thm:LichSharpChar}, see Theorem~\ref{thm:LichSharpChar} and \cite[Theorem~6.5]{cushing2018rigidity}.
We now prove, that one can drop the additional spectral condition. We specifically prove effective Bonnet Myers sharpness which implies Lichnerowicz sharpness by Theorem~\ref{thm:BMsharpImpliesReflective} and  Theorem~\ref{thm:ReflectiveImpliesLichSharp}.
\begin{theorem}\label{thm:LichSharpDistRegImpliesEffBonnMyersSharp}
Let $G=(V,E)$ be Lichnerowicz sharp and distance regular. Then, $G$ is effective Bonnet Myers sharp.
\end{theorem}

\begin{proof}
Let $f \in \R^V$ be an eigenfunction to eigenvalue $K:=\min_{x\sim y} \kappa(x,y)$. Assume $f$ is minimal in $x \in V$. We recall that spheres of radius $n$ are denoted by $S_n$.
We define another function $g \in \R^V$ via
\[
g(y) := \frac 1{|S_{d(x,y)}(x)|} \sum_{z: d(x,z)=d(x,y)} f(z).
\]
By distance regularity, $\Delta g(y) = \Delta g(y')$ whenever $d(x,y)=d(x,y')$.
Moreover, $g$ also has to be an eigenfunction to eigenvalue $K$.
 By minimality of $f(x)$ and as $f$ is not constant, we see that $g$ is not constant.
We now show that $g=cd(x,\cdot) + C$ for some constants $c,C \in \R$.

W.l.o.g., $\|\nabla g\|_\infty = 1$ and $g(z)-g(y)=1$ for some $z,y \in V$. As $g$ is constant on spheres around $x$, we can assume $y \in S_n(x)$ and $z \in S_{n+1}(x)$ for some $n$, or vice versa. As the vice versa case is analogous, we restrict to the case $y \in S_n(x)$ and $z \in S_{n+1}(x)$.
We estimate
\[
K=K(g(z)-g(y))=\Delta g(y)-\Delta g(z) \geq \inf_{\substack{h(z)-h(y)=1\\ \|\nabla h\|_\infty = 1}} \Delta h(y)-\Delta h(z) =\kappa(y,z)\geq K.
\]
Choosing $h$ the maximal Lipschitz extension of $g$ on $S_{n+2}(x)$ and the minimal Lipschitz extension of $g$ on $S_{n-1}$, we see $\Delta h(z) \geq \Delta g(z)$ and $\Delta h(y) \leq \Delta g(y)$, showing that $g=h$ on $S_{n-1}(x)$ and $S_{n+2}(x)$.  
Particularly, the maximal gradient of $g$ is also attained between $S_{n-1}(x)$ and $S_n(x)$, and between $S_{n+1}(x)$ and $S_{n+2}(x)$. Using induction, it follows that $g= d(x,\cdot) + C$ for some constant $C \in \R$.

We now show that $\sum_y d(x,y) = |V|\Deg(x)/K$.
We first notice that \[-KC=-Kg(x)=\Delta g(x) = \Deg(x)\] 
showing $C=\frac {- \Deg(x)}K$.
Moreover,
\[
0=\sum \Delta g(x) = K\sum g(x) = K\sum_y d(x,y) - \Deg(x)|V|.
\]
Rearranging shows $\sum_y d(x,y) = |V|\Deg(x)/K$.
By distance regularity, this equation holds true for all base points $x$, implying $\diameff = \frac{1}{|V|^2} \sum_{x,y} d(x,y)=\Deg(x)/K = \max \Deg/K$.
This finishes the proof.
\end{proof}

\section{Classifications}\label{sec:Classifications}

We first characterize locally connected effective Bonnet Myers sharp graphs.

\begin{theorem}\label{thm:LocConnectCharEverything}
Let $G=(V,E)$ be a locally connected graph. T.f.a.e.:
\begin{enumerate}[(i)]
\item $G$ is effective Bonnet Myers sharp,
\item $G$ is reflective,
\item $G$ is distance regular and Lichnerowicz sharp,
\item $G$ is a graph of the following list:
\begin{itemize}
\item Cocktail party graphs with at least $6$ vertices,
\item Johnson graphs, 
\item Halved cubes,
\item Schläfli graph,
\item Gosset graph.
\end{itemize}
\end{enumerate}
\end{theorem}
We note that the cocktail party graph with four vertices is locally disconnected and therefore not on the list, and the cocktail party graph with two vertices is a Johnson graph, too.
\begin{proof}
Implication $(i) \Rightarrow (ii)$ follows from Theorem~\ref{thm:BMsharpImpliesReflective}.
Implication $(ii) \Rightarrow (iv)$ follows from Lemma~\ref{lem:locConnectedAndReflectiveImpliesList}. Implication $(vi) \Rightarrow (iii)$ follows from Theorem~\ref{thm:LichSharpChar}. Finally, implication $(iii) \Rightarrow (i)$ follows from Theorem~\ref{thm:LichSharpDistRegImpliesEffBonnMyersSharp}, finishing the proof.
\end{proof}

We now drop the assumption of local connectedness.
\begin{theorem}\label{thm:effBMsharpChar}
Let $G=(V,E)$ be a graph. T.f.a.e.:
\begin{enumerate}[(i)]
\item $G$ is effective Bonnet Myers sharp,
\item $G$ is reflective and has constant curvature,
\item $G$ is a graph of the following list:
\begin{itemize}
\item Cocktail party graphs,
\item Johnson graphs, 
\item Halved cubes,
\item Schläfli graph,
\item Gosset graph,
\item Cartesian product of  graphs above with same curvature.
\end{itemize}
\end{enumerate}
\end{theorem}
\begin{proof}
Implication $(i) \Rightarrow (ii)$ follows from Theorem~\ref{thm:BMsharpImpliesReflective} and Lemma~\ref{lem:BMsharpDistanceEigenfunction}.
Implication $(ii) \Rightarrow (iii)$ follows from Theorem~\ref{thm:ReflectiveImpliesList}. We finally prove $(iii) \Rightarrow (i)$. By Theorem~\ref{thm:LocConnectCharEverything}, we have $G=G_1 \times \ldots \times G_n$, and all factors are effective Bonnet Myers sharp. As $\diameff(G)=\sum_{i=1}^n\diameff(G_i)$ and $\Deg_G = \sum_{i=1}^n \Deg_{G_{i}}$, and as all factors have the same curvature, we conclude that $G$ is also effective Bonnet Myers sharp finishing the proof.
\end{proof}

We finally drop the constant curvature assumption and characterize reflective graphs.
\begin{theorem}
Let $G=(V,E)$ be a graph. T.f.a.e.:
\begin{enumerate}[(i)]
\item $G$ is reflective,
\item $G$ is a graph of the following list:
\begin{itemize}
\item Cocktail party graphs,
\item Johnson graphs, 
\item Halved cubes,
\item Schläfli graph,
\item Gosset graph,
\item Cartesian product of graphs above.
\end{itemize}
\end{enumerate}
\end{theorem}
\begin{proof}
Implication $(i) \Rightarrow (ii)$ is given in Theorem~\ref{thm:ReflectiveImpliesList}. We finally prove $(ii) \Rightarrow (i)$. By Theorem~\ref{thm:LocConnectCharEverything}, we have $G=G_1\times\ldots \times G_n$, and all factors are reflective. By Lemma~\ref{lem:CartesianProductInheritsReflective}, this implies that $G$ is reflective finishing the proof. 
\end{proof}

\section*{Acknowledgments}
The author wants to thank Justin Salez for pointing out the relevance of the effective diameter. He moreover wants to thank Supanat Kamtue, Jack Koolen, Shiping Liu, and Norbert Peyerimhoff for useful and inspiring discussions.

\printbibliography

\appendix
\section{Bakry Emery curvature and effective diameter}

Instead of a lower bound to the Ollivier curvature, we now assume a lower bound to the Bakry Emery curvature to establish an upper bound to the effective diameter. This question however turns out to bring not many new insights as the graphs for which the effective diameter bound is attained, are precisely the hypercubes, and this is an almost immediate consequence of \cite{liu2017rigidity} as we demonstrate below. But first, we give the definition of Bakry Emery curvature (see \cite{schmuckenschlager1998curvature,lin2010ricci}).
We inductively define $\Gamma_i:\R^V\times \R^V \to \R^V$ via $\Gamma_0(f,g):=fg$ and
\[
2\Gamma_{i+1}(f,g) := \Delta \Gamma_i(f,g) - \Gamma_i(f,\Delta g) - \Gamma_i(\Delta f, g). 
\]
We write $\Gamma := \Gamma_1$ and $\Gamma_i f := \Gamma_i(f,f)$.
The Bakry Emery curvature $K(x)$ of a  vertex $x$ is given by
\[
K(x):= \inf_{\Gamma f(x) =1} \Gamma_2 f(x).
\] 
We now give the effective diameter bound.

\begin{theorem}\label{thm:BakryEmeryDiameterBound}
Let $G=(V,E)$ be a graph with Bakry Emery curvature at least $K>0$. Then,
\[
\diameff(G) \leq \frac{\max \Deg}K
\]
\end{theorem}
The proof is similar to the proof of the diameter bounds in \cite{liu2016bakry}.
\begin{proof}
We aim to show that for all $x \in V$,
\[
\frac 1 {|V|} \sum_{y\sim x} d(x,y) \leq \frac{\max \Deg}K.
\]
Let $f:=\Deg(x)$. By the gradient estimate from \cite{lin2015equivalent},
\[
\frac{(\Delta P_t f)^2}{2\Deg} \leq \Gamma P_t f \leq e^{-2Kt} P_t \Gamma f \leq e^{-2Kt} \frac{\max \Deg}2
\]
and thus,
\[
\Delta P_t f(x) \leq e^{-Kt} \max \Deg.
\]
As in the proof of Theorem~\ref{thm:EffDiamBound}, this implies 
the first claim in the proof,
and averaging over all $x \in V$ proves the theorem.
\end{proof}

\begin{theorem}
Let $G=(V,E)$ be a graph with Bakry Emery curvature at least $K>0$. Assume
\[
\diameff=\frac{\max \Deg} K.
\]
Then, $G$ is a hypercube.
\end{theorem}

\begin{proof}
Let $f=d(x,\cdot)$ for some $x\in V$. 
Due to the proof of Theorem~\ref{thm:BakryEmeryDiameterBound}, we have
\[
e^{-2Kt}P_t \Gamma f(x) = \Gamma P_t f(x)
\]
for all $t\geq  0$. From the proof of \cite[Theorem~3.4]{liu2017rigidity}, it follows that the above equality holds everywhere. Then, \cite[Theorem~4.1]{liu2017rigidity} implies that $\diam(G)=\frac {2\max \Deg} K$. Finally, \cite[Theorem~1.4]{liu2017rigidity} shows that $G$ is a hypercube  finishing the proof.
\end{proof}

\end{document}